\documentclass[12pt]{article}
\usepackage{amsfonts}
\newcommand{\oper}[2]{\newcommand{#1}{\mathop{\mathrm{#2}}\nolimits} }
\oper{\tr}{tr} \oper{\adj}{adj} \oper{\Div}{div} \oper{\ad}{ad}
\oper{\Ad}{Ad} \oper{\End}{End} \oper{\Hom}{Hob} \oper{\Aut}{Aut}
\oper{\SO}{SO} \oper{\SP}{Sp} \oper{\SU}{SU} \oper{\GL}{GL}
\oper{\T}{T} \oper{\U}{U} \oper{\id}{I} \oper{\ext}{Ext}
\oper{\rank}{rank} \oper{\diag}{Diag} \oper{\image}{image}
\oper{\zero}{Zero}

\def\cb{{\cal B}}
\def\cd{{\cal D}}

\def\ck{{\cal K}}

\def\co{{\cal O}}

\def\ct{{\cal T}}

\def\dbar{\overline\partial}

\newcommand{\CC}{\mathbb{C}}

\newcommand{\ZZ}{\mathbb{Z}}
\newcommand{\NN}{\mathbb{N}}

\def\oomega{\overline\omega}
\def\orho{\overline\rho}

\def\Oomega{\overline\Omega}

\newcommand{\lie}[1]{\mathfrak{#1}}
\newcommand\bra[2]{[{#1}\bullet {#2}]}
\newcommand\brah[2]{[{#1}\bullet {#2}]_{\lie{h}}}

\newtheorem{theorem}{Theorem}

\newtheorem{proposition}[theorem]{Proposition}
\newtheorem{definition}[theorem]{Definition}
\newtheorem{lemma}[theorem]{Lemma}

\newcommand{\bproof}{\noindent{\it Proof: }}
\newcommand{\eproof}{\  q.~e.~d. \vspace{0.2in}}
\begin{document}

\title{Extended Deformation of Kodaira Surfaces}
\author{ Yat Sun Poon \thanks{
    Partially supported by NSF DMS-0204002.}  }
\maketitle

\abstract{We present the extended Kuranishi space for Kodaira
surface as a non-trivial example to Kontsevich and Barannikov's
extended deformation theory. We provide a non-trivial example of
Hertling-Manin's weak Frobenius manifold. In addition, we find
that Kodaira surface is its own mirror image in the sense of
Merkulov. The calculations of extended deformation and the weak
Frobenius structure are based on Merkulov's perturbation method.
Our computation of cohomology is done in the context of compact
nilmanifolds.}

\section{Introduction}
Since the emergence of mirror symmetry in theoretical physics,
there are two particularly important mathematical approaches to
it. One is Kontsevich's theory of homological mirror symmetry
including a theory of extended deformation \cite{BK} \cite{Fukaya}
\cite{Kon}. Another is
 the SYZ conjecture \cite{SYZ}.

In the homological approach, the extended deformation theory
involves the full cohomology ring of complex manifolds with
coefficients in polyvector fields. Its development includes a
consideration of Frobenius manifolds and weak Frobenius manifolds
\cite{Manin2}. Although the chief impetus to both the homological
mirror symmetry and the SYZ conjecture is a desire to understand
the mathematical rationales and validity of mirror symmetry among
Calabi-Yau manifolds, the scope of Barannikov-Kontsevich theory
could be significantly extended to include complex manifolds and
symplectic manifolds in general. Work in such spirit could be
found among Merkulov's recent work
 \cite{Mer-functor} \cite{Mer-semi} \cite{Mer-a-note}
\cite{Mer-oper} as well as in Cao and Zhou's work \cite{Cao}. In
particular, Merkulov develops a homotopy version of
Hertling-Manin's weak Frobenius manifold theory and discovers an
$F_\infty$ functor. In \cite{Mer-Fro}, Merkulov considers two
objects mirror image of each other if their images through the
$F_\infty$ functor agree. The key in this analysis is concerned
with the deformation of the differential Gerstenhaber algebra
controlling the deformations of various geometric objects. We
shall briefly review the definition of differential Gerstenhaber
algebra in the next section.

In the development of extended deformation, very few examples are
known \cite{Kon}. Perhaps, this is not surprising because even
within classical deformation theory, analysis of explicit examples
of moduli space of complex manifolds has never been
straightforward. In this paper,  we present the extended Kuranishi
space of Kodaira surfaces as a stratified superspace, and find its
mirror symplectic structure in the sense of Merkulov. We shall
present Kodaira surfaces as compact quotients of a nilpotent Lie
group with a left-invariant complex structure. Geometrically, the
Kodaira surfaces are realized as an elliptic fibration over
elliptic curves. Such description allows a generalization to
higher dimensions \cite{GMPP}.

\

\noindent{\bf Main Theorem I.\ } The extended Kuranishi space of a
primary Kodaira surface
 has four components. There are  two linear components
$\CC^{6|4}$ and $\CC^{5|4}$. It has two additional components
contained in the complement of an odd linear hyperplane in
$\CC^{5|5}$ and in $\CC^{4|4}$ respectively. It has a non-trivial
weak Frobenius structure.

\

The computation of Kuranishi space is completed in Theorem
\ref{kur components}. A description of the weak Frobenius
structure is found in Section \ref{find mu}. In the above theorem,
we use the notion of supermanifold \cite{Manin1}. Our second main
result is concerned with mirror symmetry of Kodaira surfaces.

\

\noindent{\bf Main Theorem II.\ } Let $J$ be the complex structure
of a primary Kodaira surface $N$. There exists a family of
symplectic structures on the underlying smooth manifold of $N$
such that each symplectic structure in this family is the mirror
image of $(N, J)$.

\

The proof of Theorem II  is based on an analysis of the
differential Gerstenhaber algebras controlling the deformations of
complex structures and symplectic structures. Its proof is
completed in Theorem \ref{follow-up}. The key ingredients are
Lemma \ref{complex dG} and Lemma \ref{symp dG} as they allow us to
translate the issue to a finite-dimensional algebraic problem.

\section{Differential Gerstenhaber Algebra of Nilmanifolds}

 In this
section, we recall the definition and elementary facts of
differential Gerstenhaber algebras, and construct a collection of
examples based on 2-step nilpotent algebra with a complex
structure.

\subsection{Differential Gerstenhaber Algebra (DGA)}\label{def of G}

Let $\lie{f}=\oplus_{j\in\NN }\lie{f}^j$ be a complex graded
vector space. If $v$ is an element in $\lie{f}^j$, it is said to
be homogeneous and we denote its grading $j$ by $\tilde v$.

\begin{definition} A Gerstenhaber algebra is a graded vector space
$\lie{f}$ equipped with two product structures, a Schouten bracket
$\bra{}{}$ and a wedge product $\wedge$, such that the following
set of axioms hold.
\begin{itemize}
\item[{\rm (L1)}] $\bra{\lie{f}^i}{\lie{f}^j}\subset
\lie{f}^{i+j-1}$, \item[{\rm (L2)}]
$\bra{a}{b}=(-1)^{\tilde{a}\tilde{b}+\tilde{a}+\tilde{b}}\bra{b}{a}$,
 \item[{\rm (L3)}]
$\bra{a}{\bra{b}{c}}=\bra{\bra{a}{b}}{c} -(-1)^{{\tilde a}{\tilde
b}+\tilde{a}+\tilde{b}} \bra{b}{\bra{a}{c}}$, \item[{\rm (C1)}]
$\lie{f}^i\wedge\lie{f}^j\subset \lie{f}^{i+j}$,
 \item[{\rm (C2)}] $a\wedge b
=(-1)^{{\tilde a}{\tilde b}}b\wedge a$, \item[{\rm (C3)}]
$\bra{a\wedge b}{c}=
a\wedge\bra{b}{c}+(-1)^{\tilde{a}\tilde{b}}b\wedge\bra{a}{c}$.
\end{itemize}
 If in addition, there is an operator $\dbar$ such that
\begin{itemize}
\item[{\rm (D1)}] $\dbar \lie{f}^j\subset \lie{f}^{j+1}$,
\item[{\rm (D2)}] $\dbar\circ \dbar=0$, \item[{\rm (D3)}] $\dbar
\bra{a}{b}=\bra{\dbar a}{b}-(-1)^{\tilde a}\bra{a}{\dbar b}$,
\item[{\rm (D4)}] $\dbar(a\wedge b) =(\dbar a)\wedge
b+(-1)^{{\tilde a}}a\wedge (\dbar b)$,
\end{itemize}
then the collection $(\lie{f}, \bra{}{}, \wedge, \dbar)$ is said
to be a differential Gerstenhaber algebra, or simply DGA.
\end{definition}

If we ignore the distributive law (C3), the structure of a DGA
could be considered as two entities. The triple $(\lie{f}, \wedge,
\dbar)$ with axioms (C1), (C2), (D1), (D2) and (D4) forms a
differential graded algebra. The triple $(\lie{f}, \bra{}{},
\dbar)$ with axioms (L1), (L2), (L3), (D1), (D2) and (D3) forms a
differential graded Lie algebra. If the grading of the algebra
$\lie{f}$ is reduced to $\ZZ_2$, it is the structure of odd
differential Lie superalgebra \cite[page 155]{Manin1} \cite[page
160]{Manin2}. (L3) should be considered as the Jacobi identity.
Its alternative expression is
\[
(-1)^{(\tilde{a}+1)(\tilde{c}+1)}\bra{a}{\bra{b}{c}}+
(-1)^{(\tilde{b}+1)(\tilde{a}+1)}\bra{b}{\bra{c}{a}}+
(-1)^{(\tilde{c}+1)(\tilde{b}+1)}\bra{c}{\bra{a}{b}}=0.
\]
The commutative law (L2) and the distributive law (C3) could also
be expanded to include the following.
\begin{equation}
\bra{a}{b}=(-1)^{\tilde{a}\tilde{b}+\tilde{a}+\tilde{b}}\bra{b}{a}=
 -(-1)^{(\tilde{a}+1)(\tilde{b}+1)}\bra{b}{a}.\label{commutative}
 \end{equation}
 \begin{equation}
\bra{a}{b\wedge c}=\bra{a}{b}\wedge
c+(-1)^{\tilde{b}\tilde{c}}\bra{a}{c}\wedge b= \bra{a}{b}\wedge
c+(-1)^{\tilde{b}+\tilde{a}\tilde{b}}b\wedge\bra{a}{c} .
\label{distributive1}
\end{equation}
\begin{equation}
 \bra{a\wedge
b}{c}=a\wedge\bra{b}{c}+(-1)^{\tilde{a}\tilde{b}}b\wedge\bra{a}{c}=
a\wedge\bra{b}{c}+(-1)^{\tilde{b}+\tilde{b}\tilde{c}}\bra{a}{c}\wedge
b . \label{distributive2}
\end{equation}

Given a differential Gerstenhaber algebra, we may consider the
cohomology with respect to the operator $\dbar$. Due to axioms
(D3) and (D4), the Schouten bracket and the wedge product descend
to the cohomology of the complex so that the cohomology becomes a
Gerstenhaber algebra.

\begin{definition} Suppose that $dG$ and $\widehat{dG}$ are two
differential Gerstenhaber algebras. A homomorphism $\Upsilon$ from
$dG$ to $\widehat{dG}$ is a quasi-isomorphism if the induced map
from the cohomology of $dG$ to the cohomology of $\widehat{dG}$ is
an isomorphism of Gerstenhaber algebras.
\end{definition}

\subsection{Gerstenhaber algebra associated to a Lie
algebra with complex structures}\label{construction of G}

 Suppose that $\lie{g}$ is a real finite-dimensional Lie algebra with Lie bracket $[\ , \ ]$. Let $J$ be an
integrable complex structure on $\lie{g}$. In other words, $J\circ
J=-\mbox{identity}$,  and for any pair of elements $A$ and $B$ in
$\lie{g}$
\begin{equation}
[A, B]-[JA, JB]+J[JA, B]+J[A, JB]=0.
\end{equation}

The complexification of $\lie{g}$ is denoted by $\lie{g}_\CC$. It
decomposes into type (1,0)-vectors and (0,1)-vectors. i.e. $
\lie{g}_\CC=\lie{g}^{1,0}\oplus \lie{g}^{0,1}$.  The integrability
assumption implies that the vector subspaces $\lie{g}^{1,0}$ and
$\lie{g}^{0,1}$ are subalgebras. If $\omega$ is an element in
$\wedge^k\lie{g}^*$, define
\begin{equation}
(J\omega)(A_1, \dots, A_k)=(-1)^k\omega(JA_1, \dots, JA_k).
\end{equation}
We obtain a complex structure on the dual vector space
$\lie{g}^*$. It has a corresponding decomposition $
\lie{g}^*_\CC=\lie{g}^{*(1,0)}\oplus\lie{g}^{*(0,1)}. $ Let
$(\lie{f}, \wedge)$ be the exterior algebra generated by the
vector space $\lie{g}^{*(0,1)}\oplus\lie{g}^{1,0}$. The degree-k
elements form a vector subspace
$\lie{f}^k:=\oplus_{p+q=k}\lie{g}^{*(0,q)}\otimes\lie{g}^{p,0}$
where $\lie{g}^{p,0}$ is the $p$-th exterior product of
$\lie{g}^{1,0}$. The space $\lie{g}^{*(0,q)}$ is the $q$-th
exterior product of its first degree counterpart.

Next, we generate the Schouten bracket $\bra{}{}$ on $\lie{f}$ in
the following way. If $U_j$ and $V_\ell$ are in $\lie{g}^{1,0}$,
then $\bra{U_j}{V_\ell}=[U_j, V_\ell]$ is the Lie bracket of the
algebra $\lie{g}$. In general,
\begin{eqnarray*}
 & &\bra{U_1\wedge\cdots\wedge U_k}{V_1\wedge\cdots\wedge V_p}\nonumber\\
&:=&\sum_{j=1}^{k}(-1)^{(k-j)} U_1\wedge\cdots \wedge
U_{j-1}\wedge U_{j+1}\wedge \cdots\wedge U_{k} \wedge \bra{U_j}{
V_1\wedge \cdots\wedge V_p}.
\end{eqnarray*}
For $\psi$ in $\lie{g}^{*(0,p)}$ and $V_1\wedge\cdots\wedge V_k$
in $\lie{g}^{k,0}$, define
\[
\bra{ V_1\wedge\cdots\wedge V_k}{\psi}:=\sum_{j=1}^k(-1)^{k-j}
V_1\wedge \cdots\wedge V_{j-1}\wedge V_{j+1}\wedge\cdots\wedge
V_k\wedge (L_{V_j}\psi).
\]
If $\phi$ is in $\lie{g}^{*(0,\ell)}$, we define
$
\bra{\phi}{\psi}=0.
$
Finally for $\phi\wedge \Xi$ in
$\lie{g}^{*(0,q)}\otimes\lie{g}^{p,0}$ and $\psi\wedge \Theta$ in
$\lie{g}^{*(0,\ell)}\otimes\lie{g}^{k,0}$, we define
\begin{eqnarray}
&&\bra{\phi\wedge \Xi}{\psi\wedge \Theta}\nonumber\\
&:=&\phi\wedge\bra{\Xi}{\psi}\wedge\Theta
     +(-1)^{(\tilde\phi+\tilde\Xi)(\tilde\psi+\tilde\Theta)+
(\tilde\phi+\tilde\Xi)+(\tilde\psi+\tilde\Theta)}
    \psi\wedge\bra{\Theta}{\phi}\wedge\Xi\nonumber\\
&&+(-1)^{\tilde\psi(\tilde\Xi+1)}\phi\wedge\psi\wedge\bra{\Xi}{\Theta}.
\end{eqnarray}

It is straightforward to verify that the structure $(\lie{f},
\bra{}{}, \wedge)$ forms a Gerstenhaber algebra.

\subsection{Construction of a differential
operator}\label{construction of dG}

Suppose that the complex structure $J$ is abelian. In other words,
 for any $A$ and $B$ in $\lie{g}$,
\begin{equation}\label{abelian}
[JA, JB]=[A, B].
\end{equation}
Such complex structure is characterized by the eigenspaces
$\lie{g}^{1,0}$ and $\lie{g}^{0,1}$ with respect to $J$ being
abelian subalgebras of the complexified algebra $\lie{g}_\CC$
\cite{Dotti}.

 Suppose in addition that the Lie algebra $\lie{g}$ is a 2-step nilpotent Lie
algebra. The quotient $\lie{t}$ of the algebra $\lie{g}$ by the
center $\lie{c}$ is abelian.
 It yields an exact sequence of algebras
\begin{equation}\label{central extension}
0\to \lie{c}\to \lie{g}\to \lie{t} \to 0.
\end{equation}
As vector spaces, there is a  direct sum decomposition
$\lie{g}_\CC=\lie{t}_\CC\oplus\lie{c}_\CC$. Since the complex
structure $J$ is abelian, the center is invariant with respect to
$J$. Therefore, as vector spaces,
\begin{equation}
\lie{g}^{1,0}=\lie{t}^{1,0}\oplus\lie{c}^{0,1}, \quad
\lie{g}^{0,1}=\lie{t}^{0,1}\oplus\lie{c}^{0,1}.
\end{equation}
Let $\{T_j, 1\leq j\leq n\}$ be a basis for $\lie{t}^{1,0}$.

Assume that the real dimension of the center $\lie{c}$ is
two-dimensional. Let $W$ be a non-zero element in $\lie{c}^{1,0}$.
Then the 2-step nilpotence and (\ref{abelian}) together imply that
the structural equations of the Lie algebra are given by some
constants $E_{kj}$ and $F_{kj}$ such that
\begin{equation}\label{structural}
[\bar{T}_j, {T}_k]= E_{jk} W+ F_{jk} \bar{W}.
\end{equation}
The structural constants are subjected to the constraint
$\bar{F}_{kj}=-E_{jk}$.

Let $\{\omega^j, 1\leq j\leq n\}$ and $\{\rho\}$ be the dual bases
of $\{T_j\}$ and $\{W\}$ respectively. Then the dual structural
equations are
\begin{equation}
d\omega^j=0, \quad
d\rho=-\sum_{j,k}E_{jk}\oomega^j\wedge\omega^k=\sum_{j,k}E_{jk}\omega^k\wedge\oomega^j.
\end{equation}
It follows that
\begin{equation}\label{dual structural}
\dbar\oomega^j=0, \quad \dbar\orho=0, \quad \mbox{ and } \quad
\dbar\rho=\sum_{j,k}E_{jk}\omega^k\wedge\oomega^j.
\end{equation}

 We define a linear map $\dbar: \lie{g}^{1,0}\to
\lie{g}^{*(0,1)}\otimes\lie{g}^{1,0}$ as follows \cite{Gauduchon}.
For any $(1,0)$-vector $A$ and $(0,1)$-vector $\bar{B}$, define
\begin{equation}
\dbar_{\bar{B}}A:=[{\bar{B}}, A]^{1,0}.
\end{equation}
It follows that $\dbar W=0$ and $\dbar T_j=E_{kj}\oomega^k\wedge
W.$   Using the usual $\dbar$-operator on differential forms and
extension by anti-derivation on polyvectors, we  extend this
operator to $\lie{f}$. To be precise, for $\Oomega\wedge \Xi$ in
$\lie{g}^{*(0,q)}\otimes\lie{g}^{p,0}$,
\begin{equation}
\dbar (\Oomega\wedge \Xi)=\dbar \Oomega \wedge
\Xi+(-1)^q\Oomega\wedge\dbar \Xi.
\end{equation}
Given the $\overline{\partial}$-operator, we have a resolution for
each $\lie{g}^{p,0}$.
\begin{equation}\label{complex}
0\to \lie{g}^{p,0} {\rightarrow}
\lie{g}^{*(0,1)}\otimes\lie{g}^{p,0} {\rightarrow} \cdots
\lie{g}^{*(0,k)}\otimes\lie{g}^{p,0} {\rightarrow}
\lie{g}^{*(0,k+1)}\otimes\lie{g}^{p,0} {\rightarrow} \cdots.
\end{equation}
\begin{lemma}When the dimension of $\lie{c}^{1,0}$ is equal to one,
the above resolution is a complex.
\end{lemma}
\bproof Since $\dbar\Oomega=0$ for all $\Oomega$ in
$\lie{g}^{*(0,q)}$,  it suffices to prove that $\dbar^2 \Xi=0$ for
all polyvector field $\Xi$. Since $\dbar W=0$, we focus on $\Xi
\in \lie{t}^{k,0}$. By re-ordering elements in the basis, it
suffices to consider $\Xi=T_1\wedge \cdots \wedge T_k$.
\begin{eqnarray}
& &\dbar\dbar(T_1\wedge \cdots \wedge T_k)\nonumber\\
&=&\sum_{j=1}^k E_{kj}
 \dbar(\oomega^k\wedge T_1\wedge \cdots \wedge W \wedge \cdots \wedge T_k )
\nonumber\\
&=&-\sum_{j=1}^k E_{kj} \oomega^k\wedge\dbar(T_1\wedge \cdots
\wedge W \wedge \cdots \wedge T_k).
\end{eqnarray}
Since the vector part of $\dbar T_j$ is the linear span of $W$ due
to dimension restriction and $\dbar W=0$ due to the step of
nilpotence, the last operation yields zero. \eproof

\begin{definition}\label{cohomology} Let $\lie{h}^{p, q}$ be the $q$-th cohomology of the complex
{\rm (\ref{complex})}. i.e.
\begin{eqnarray*}
\lie{h}^{p, q} &=&\frac{ \ker \dbar:
\lie{g}^{*(0,q)}\otimes\lie{g}^{p,0} {\rightarrow}
\lie{g}^{*(0,q+1)}\otimes\lie{g}^{p,0}} { \dbar \left(
\lie{g}^{*(0,q-1)}\otimes\lie{g}^{p,0}\right)}
=\frac{\ker\dbar_q}{\mbox{\rm \rm image }\dbar_{q-1}}.\\
\lie{h}^k&=&\oplus_{p+q=k}\lie{h}^{p,q}.
\end{eqnarray*}
\end{definition}

It is now straightforward to verify the following.

\begin{proposition}\label{existence of dGA}
For any 2-step nilpotent algebra $\lie{g}$ with an abelian complex
structure $J$ and real two-dimensional center, $(\lie{f},
\bra{}{}, \wedge, \dbar)$ is a differential Gerstenhaber algebra.
\end{proposition}

Concerned with the scope of the last theorem, we note that there
are at least six series of Lie algebras qualified for the
assumptions \cite{Salamon-nil}.

\subsection{2-Step nilmanifolds}\label{nil}

Let $\lie{g}$ be a Lie algebra as described in Proposition
\ref{existence of dGA}. Let $G$ be a simply connected Lie group
for the algebra $\lie{g}$. Suppose that there exists a discrete
subgroup $\Gamma$ of $G$ such that the quotient space
$N:=\Gamma\backslash G$ with respect to the left action of
$\Gamma$ is compact.  The resulting quotient is a nilmanifold. On
the group $G$, we use left translation of $J$ to define a
left-invariant almost complex structure on the manifold $G$. Due
to invariance, it descends to a complex structure on the quotient
space $N$. Since the almost complex structure is left invariant
and abelian, its Nijenhuis tensor vanishes and $N$ becomes a
compact complex manifold.

The extension (\ref{central extension}) yields  a principal
holomorphic fibration of $N$ over an abelian variety with elliptic
curves as fibers.
\begin{equation}\label{fibration}
\Psi: N=\Gamma\backslash G \to M=T^{2n}.
\end{equation}
The vector field generated by $W$ is a global holomorphic vector
field trivializing the kernel of the differential of the
projection map $\Psi$. Let $\ct_N$ and $\co_N$ be the holomorphic
vector bundle and the structure sheaf respectively for the complex
manifold $N$. The projection above yields an exact sequence of
holomorphic vector bundles on $N$.
\begin{equation}\label{extension-1}
0\to \lie{c}^{1,0}\otimes\co_N \to \ct_N \to \Psi^*\ct_M\to 0.
\end{equation}
Since $\lie{c}^{1,0}$ is one-dimensional, an inspection of
transition function implies that for all $p$, one has the exact
sequence
\begin{equation}\label{extension-2}
0\to \lie{c}^{1,0}\otimes\Psi^*\wedge^p\ct_M \to \wedge^{p+1}\ct_N
\to \Psi^*\wedge^{p+1}\ct_M\to 0.
\end{equation}

Based on the above exact sequence, a computation similar to those
in the proof of Theorem 1 in \cite{MPPS} shows the following.
\begin{proposition}\label{cx qua}
Suppose that $N$ is a 2-step nilmanifold with abelian complex
structure and complex one-dimensional center. The inclusion of the
invariant differential Gerstenhaber algebra $dG(\lie{g}, J)$ in
the algebra $dG(N, J)$ is a quasi-isomorphism. I.e., there is a
natural isomorphism $H^q(N, \wedge^p\ct_N)\cong \lie{h}^{p, q}$.
\end{proposition}

To find representatives in cohomology classes, we take the
Hermitian metric on $N=\Gamma\backslash G$ such that the frame
$\{W, T_j, 1\leq j\leq n\}$ is Hermitian. This choice of Hermitian
metric induces a Hermitian metric on the bundle of polyvectors
$\wedge^p\ct_N$ and a Hermitian inner product on
$\lie{g}^{*(0,q)}\otimes\lie{g}^{p,0}$. In terms of this inner
product, we consider  the orthogonal complement of the image
$\lie{g}^{*(0,q-1)}\otimes\lie{g}^{p,0}$ through $ \dbar$ in the
kernel of $\dbar$ on $\lie{g}^{*(0,q)}\otimes\lie{g}^{p,0}$.
Denote this space by $\dbar^\perp
(\lie{g}^{*(0,q-1)}\otimes\lie{g}^{p,0})$.

\begin{theorem}\label{harmonic rep} For each element $\psi\wedge \Xi$ in
$\lie{g}^{*(0,q)}\otimes \lie{g}^{p,0}$, $\dbar^* (\psi\wedge
\Xi)$ with respect to the $L_2$-norm on the compact Hermitian
manifold $N$ is equal to $\dbar^*(\psi\wedge\Xi)$ with respect to
the Hermitian inner product on the finite-dimensional vector
spaces $\lie{g}^{*(0,q)}\otimes \lie{g}^{p,0}$. In particular, $
\dbar^\perp (\lie{g}^{*(0,q-1)}\otimes\lie{g}^{p,0}) $ is a space
of harmonic representatives for Dolbeault cohomology $H^q(N,
\wedge^p\ct_N)$ on the compact complex manifold $N$.
\end{theorem}

The idea behind the proof of this theorem is already in
\cite{GMPP} and \cite{MPPS}. We do not repeated the details here.

\subsection{Kodaira manifolds}\label{K-man}
Kodaira surfaces are locally trivial elliptic fibrations over
elliptic curves \cite{Kodaira}. They could be realized as the
compact quotient of a nilpotent extension of a three-dimensional
Heisenberg group. Below, we extend this description to all
dimensions \cite{GMPP}.

Consider the manifold $R^{2n+2}$ with coordinates $(x_j, y_j, u,
v)$, where $1\leq j\leq n$. Define a multiplication by
\begin{eqnarray}
& &(x_j,y_j, u, v)\ast (x_j^{\prime },y_j^{\prime },u^{\prime },
v^\prime)\nonumber\\
&=&\left(x_j+x_j^{\prime},
    y_j+y_j^{\prime},u+u^{\prime}+
    \frac{1}{2}\sum_{j=1}^n(x_jy_j^{\prime}-y_jx_j^{\prime}), v+v^\prime
\right).
\end{eqnarray}
This multiplication turns $R^{2n+2}$ into a Lie group with the
origin as identity element. This group is the product of the
(2n+1)-dimensional Heisenberg group and the one-dimensional
additive abelian group. The former is given by $v=0$. The latter
is given by $x_j=y_j=u=0$ for all $j$. Let $(X_j, Y_j, U, V)$ be
the left invariant vector fields generated by left translation of
$\frac{\partial}{\partial x_j }$, $\frac{\partial }{\partial
y_j}$, $\frac{\partial}{\partial u}$, $\frac{\partial}{\partial
v}$ at the identity element. Its algebra is denoted by $\lie{g}$.
The non-zero structural equations are given by $ [X_j,Y_j]=U$ for
all $1\leq j\leq n$. The algebra generated by $X_j, Y_j$ and $U$
is the Heisenberg algebra. The center $\lie{c}$ of $\lie{g}$ is
spanned by $U$ and $V$. The quotient of $\lie{g}$ with respect to
its center is denoted by $\lie{t}_{2n}$.

On the algebra $\lie{g}$, define an endomorphism $J$ by
\begin{equation}\label{defining J}
JX_j=Y_j, \quad JY_j=-X_j, \quad JU=V, \quad JV=-U.
\end{equation}
It is an abelian complex structure. Let $ T_j=\frac12(X_j-iY_j)$
and $W=\frac12(U-iV).$
 Then the non-zero
structural constants are $ E_{jj}=-\frac{i}2$ for $1\leq j\leq n$.
Based on this complex structure, we construct the corresponding
differential Gerstenhaber algebra $\lie{f}$ as in the previous
sections. The following two observations will be useful for our
future computation.
\begin{lemma}\label{2-dbar}
The $\dbar$-operator of the Gerstenhaber algebra $\lie{f}$ is
completely generated by $ \dbar W=0$, $\dbar
T_j=-\frac{i}2\oomega^j\wedge W$, $\dbar\oomega^j=0$, and
$\dbar\orho=0$.
\end{lemma}
\begin{lemma}\label{t-rho}
The Schouten bracket in the Gerstenhaber algebra $\lie{f}$ is
completely generated by identity
$\bra{T_j}{\orho}=-\frac{i}2\oomega^j$.
\end{lemma}

Let $e_1, \dots, e_{2n}, e_{2n+1}, e_{2n+2}$ be the standard basis
for the vector space $R^{2n+2}$. The rank of the discrete subgroup
$\Gamma$ generated by them is equal to 2n+1. Its quotient is a
compact manifold. This is a Kodaira manifold. We denote it by $N$.
Given Lemma \ref{2-dbar} and Lemma \ref{t-rho}, we  may calculate
cohomology spaces. For example, the cohomology spaces below are
relevant to the deformation of generalized complex structures
\cite{Gua}.
\begin{lemma}\label{degree-2} On a Kodaira manifold,
the degree-two harmonic fields in $\lie{h}$ are given as follows.
\begin{eqnarray*}
\lie{h}^{0,2}&=&\langle \cb^{j}=\oomega^j\wedge\orho, \quad \cb^{ij}=\oomega^i\wedge\oomega^j\rangle,
\quad \lie{h}^{2,0} =\langle B_j=T_j\wedge W \rangle,\\
 \lie{h}^{1,1}&=&\langle \phi=\orho\wedge W, \quad
\phi^j_k=\frac12(\oomega^j\wedge T_k+\oomega^k\wedge T_j)\rangle.
 \end{eqnarray*}
\end{lemma}

We could also identify some elements in the center with respect to
the Schouten bracket. Note that $d\oomega_j=0$ for all $1\leq
j\leq n$, and $W$ commutes with all vectors, it follows that
$\oomega_j$ and $W$ commute with all elements in the algebra
$\lie{f}$ with respect to the Schouten bracket. Furthermore, let
us consider elements of the form
$\oomega^1\wedge\cdots\wedge\oomega^n\wedge \phi\wedge\Xi$ where
$\phi\wedge\Xi\in \lie{g}^{*(0,\ell)}\otimes\lie{g}^{k,0}$. For
any $\psi\wedge \Theta$ in $\lie{g}^{*(0,q)}\otimes\lie{g}^{p,0}$,
\[
\bra{\oomega^1\wedge\cdots\wedge\oomega^n\wedge
\phi\wedge\Xi}{\psi\wedge\Theta}
=\oomega^1\wedge\cdots\wedge\oomega^n\wedge \bra{
\phi\wedge\Xi}{\psi\wedge\Theta}
\]
Since $\bra{\Xi}{\psi}$ is contained in
$\lie{t}^{*(0,q)}\otimes\lie{g}^{k-1,0}$, $\bra{\Theta}{\phi}$ is
contained in $\lie{t}^{*(0, \ell)}\otimes\lie{g}^{(p-1,0)}$, and
$\oomega^1\wedge\cdots\wedge\oomega^n$ spans the highest exterior
product of $\lie{t}^{*(0,1)}$,
 the right hand side of the above equality is equal to zero so long as $\ell>0$ and
 $q>0$.
\begin{lemma}\label{full rank} Elements of the form
$\oomega^1\wedge\cdots\wedge\oomega^n\wedge \phi\wedge\Xi$ where
$\phi\wedge\Xi\in \lie{g}^{*(0,\ell)}\otimes\lie{g}^{k,0}$ with
$\ell>0$ are in the center of the algebra $(\lie{f},\bra{}{})$.
The subspaces $\lie{t}^{*(0, \ell)}$ and $\lie{t}^{*(0,
q)}\otimes\lie{c}^{1,0}$ are in the center of $(\lie{f},
\bra{}{})$.
\end{lemma}

\section{Extended Deformations of Kodaira Surfaces}

In this section, we apply the theory of extended deformations as
developed by Barannikov, Kontsevich \cite{BK} and Merkulov
\cite{Mer-Fro} to Kodaira surfaces.

\subsection{A brief overview of extended
deformation}\label{review}
 Rigorous treatment
of extended deformation theory could be found in several papers.
 Our principle reference is \cite{Mer-Fro}. To set up the notations, we present a
shortcut into the computational aspect of this theory.

On a complex manifold, there is a natural differential
Gerstenhaber algebra. Its graded vector space is the space of
smooth sections of bundle of $(0,q)$-forms with values in
$(p,0)$-vectors.
\[
\lie{n}=\oplus_k\lie{n}^k=\oplus_k\left(\oplus_{p+q=k}C^\infty(N,
\wedge^{q}T^{*(0,1)}_N\otimes \wedge^{p}T^{1,0}_N)\right).
\]
The Schouten bracket and the wedge product are fiber-wise defined
as we did for the finite-dimensional vector spaces in Section
\ref{construction of G} and Section \ref{construction of dG} .

The $\dbar$-operator is  the $\dbar$-operator on differential
forms twisted by the Chern connection of a Hermitian metric on
holomorphic tangent bundle.  The cohomology is
\begin{equation}
\lie{h}=\oplus_k\lie{h}^k, \quad \lie{h}^k=\oplus_{p+q=k}H^q (N,
\wedge^p \ct_N).
\end{equation}

In analogue with classical case, we treat elements in $\lie{h}$ as
infinitesimal deformations. To construct a deformation with
prescribed infinitesimal deformation, we take a harmonic base
$\{\theta^\alpha\}$ for $\lie{h}$. Let $\Gamma_1$ be a harmonic
representative with coordinates ${\bf x}=(x_\alpha)$. Here we
consider $\lie{h}$ as a \lq supermanifold\rq
 \ with even and odd coordinates. A coordinate $x_\alpha$ is odd if
and only if the grading of the corresponding harmonic element
$\theta^\alpha$ in the Gerstenhaber algebra $\lie{h}$ is odd. Next
we consider the formal ring $\lie{n}[[{\bf x}]]$ of power series
in $(x_\alpha)$ with coefficients in $\lie{n}$. The grading of
$\lie{n}$ together with the parity of the coordinate functions
gives elements in this ring a $\ZZ_2$-grading. For example, as an
element in $\lie{n}[[{\bf x}]]$, $\Gamma_1=\sum_\alpha
x_\alpha\theta^\alpha$ is always even. Now, $\Gamma_1$ is
integrable if there is an even element $\Gamma$ in $\lie{n}[[{\bf
x}]]$ with the following properties.
\begin{enumerate}
\item Up to degree 1, $\Gamma=\Gamma_1$. We denote it by
$\Gamma\equiv_1\Gamma_1$. \item $\Gamma$ solves the Maurer-Cartan
equation.
\end{enumerate}
Given any $\Gamma_1$, one could apply the Kuranishi recursive
method to generate a solution to the Maurer-Cartan equation and
identify the obstruction \cite{Mer-Fro}. The obstruction is stored
in an odd vector field $\vec\partial$ on the supermanifold
$\lie{h}$. Following Merkulov, we call the vector field
$\vec\partial$ the Chen vector field \cite{Chen}. This vector
field is uniquely determined by the fact that if $\Gamma$ is
generated by the Kuranishi recursive method, $\vec\partial\Gamma$
is equal to the harmonic part of $-\dbar
\Gamma-\frac12\bra{\Gamma}{\Gamma}$. Therefore, finding extended
deformation is to solve the \lq extended' Maurer-Cartan equation
\begin{equation}\label{EMC}
\dbar\Gamma+\vec\partial\Gamma+\frac12\bra{\Gamma}{\Gamma}=0
\end{equation}
with initial condition $\Gamma\equiv\Gamma_1$.

Up to gauge equivalence the extended Kuranishi space is the zeroes
of the vector field $\vec\partial$. Let $[\ ,\ ]$ be the usual
(super) Lie bracket of vector fields on the infinite-dimensional
supermanifold $\lie{n}$. The gauge equivalence is due to the
distribution
\begin{equation}
\cd:=\{[\vec\partial, Y]: Y \mbox{ is a vector field on }
\lie{n}\}.
\end{equation}
In other words, the Kuranishi space is the quotient space.
\begin{equation}
\ck:=\frac{\zero\vec\partial}{\cd}.
\end{equation}
Along the Kuranishi space, we consider the operator
$\dbar_\Gamma:=\dbar+\bra{\Gamma}{\ \ }$. It is known that
 the structure $(\lie{n}[[{\bf x}]], \bra{}{},
\wedge, \dbar_\Gamma)$ is a differential Gerstenhaber algebra
\cite[4.4]{Mer-Fro}. It means that extended deformation of complex
structures gives deformation of differential Gerstenhaber
algebras. A variation of the induced Gerstenhaber algebras on the
cohomology level follows. To describe this variation, we take a
partial derivative of the Maurer-Cartan equation with respect to
the super-coordinate $x_\alpha$.
\begin{equation}
\dbar_\Gamma\left(\frac{\partial\Gamma}{\partial
x_\alpha}\right)=\dbar\left(\frac{\partial\Gamma}{\partial
x_\alpha}\right)+\bra{\Gamma}{\frac{\partial\Gamma}{\partial
x_\alpha}}=0.
\end{equation}
Therefore, $\frac{\partial\Gamma}{\partial x_\alpha}$ represents a
$\dbar_\Gamma$-closed class. Examining the first order terms for
$\frac{\partial\Gamma}{\partial x_\alpha}$  and for all $\alpha$,
we realize that they span the entire cohomology.

As we have a deformed differential Gerstenhaber algebra,
\[
\frac{\partial\Gamma}{\partial x_\alpha}\wedge
\frac{\partial\Gamma}{\partial x_\beta}
\]
is again $\dbar_\Gamma$-closed. Therefore, up to a
$\dbar_\Gamma$-exact term there exist  functions
$\mu^{\alpha\beta}_\gamma$ such that
\begin{equation}
\frac{\partial\Gamma}{\partial x_\alpha}\wedge
\frac{\partial\Gamma}{\partial x_\beta}
=\sum_\gamma\mu^{\alpha\beta}_\gamma
\frac{\partial\Gamma}{\partial x_\gamma}
\end{equation}
Now consider a product on the tangent bundle $\ct_{\ck}$ of the
supermanifold $\ck$ defined by
\begin{equation}\label{mu-2}
\frac{\partial}{\partial x_\alpha} \circ \frac{\partial}{\partial
x_\beta}=
\sum_\gamma\mu^{\alpha\beta}_\gamma\frac{\partial}{\partial
x_\gamma}
\end{equation}
then the product varies over the superspace $\ck$. This determines
an associative and commutative product on the tangent sheaf
$\ct_\ck$ at least along the smooth part of $\ck$ \cite[Theorem
3.6.2]{Mer-Fro}. It makes $(\ct_\ck, \circ)$ a F-manifold or weak
Frobenius manifold \cite[Corollary 4.9.2]{HM}.

In the presence of obstruction, the entire computation is replaced
by the operator $ \dbar+\vec\partial+\bra{\Gamma}{\ \ }$ and the
computation is done on the supermanifold $\lie{h}$. In such case,
Merkulov finds a homotopy version of F-manifolds, so called
F$_\infty$ structure on $\lie{h}$.

\subsection{Generalized deformations of Kodaira manifolds}\label{baby}
As an example, we solve the Maurer-Cartan equation generated by
elements in $\lie{h}^2$ for Kodaira manifolds. One may consider
such deformations as deformation of generalized complex structures
\cite{Gua}.

Due to Lemma \ref{full rank}, the only non-zero brackets among
elements in $\lie{h}^2$ are
\begin{equation}
\bra{\cb^j}{B_k}= \dbar s^j_k, \mbox{ for all } 1\leq j,k\leq n,
\end{equation}
where \begin{equation} s^j_k=\frac12(\oomega^j\wedge
T_k-\oomega^k\wedge T_j) \quad \mbox{ and } \quad \dbar s^j_k =
\frac{i}2\ \oomega^j\wedge\oomega^k\wedge W. \label{ds}
\end{equation}
It follows that the field $s^j_k$ appears in the Kuranishi
recursive formula, and hence we have to calculate its bracket with
and all elements in $\lie{h}^2$. It turns out the only non-zero
brackets are the following.
\begin{equation}
\bra{s^j_k}{\cb^\ell}
=\frac{i}2\oomega^j\wedge\oomega^k\wedge\oomega^\ell, \quad
\bra{s^j_k}{\psi}= -\dbar s^j_k.
\end{equation}
Since $\oomega^j\wedge\oomega^k\wedge\oomega^\ell$ commutes with
all other elements in $\lie{f}$ with respect to the Schouten
bracket, it is now straightforward to verify the following.

\begin{theorem}\label{generalized}
Every direction $\Gamma_1$ in $\lie{h}^2$ as an infinitesimal
extended deformation is unobstructed. If
\[
\Gamma_1=\sum_{j=1}^na_j\cb^j+\sum_{1=i<j=n}a_{ij}\cb^{ij}
+a\psi+\sum_{j,k=1}^na^k_j\phi^j_k+ \sum_{j=1}^na^jB_i
\]
and $a\neq 1$, then $(\vec\partial, \Gamma)$ with
\begin{equation}
\vec\partial=\vec{0}, \quad {\mbox{ and }} \quad
\Gamma=\Gamma_1-\frac{1}{1-a}\sum_{i,j=1}^na^ia_j s^j_i
\end{equation}
form a solution to the extended Maurer-Cartan equation with
$\Gamma_1$ as initial condition. In other words, every element in
$\lie{h}^2$ for a Kodaira manifold is unobstructed.
\end{theorem}

\subsection{DGA of Kodaira surfaces}\label{full dG}

In the remaining sections, we consider the extended deformation of
 of Kodaira surfaces in all directions in $\lie{h}$.
 We begin with finding the cohomology of all
degrees.

Firstly, we consider the resolution of $\lie{g}^{1,0}$. By Lemma
(\ref{2-dbar}), $ \lie{h}^{1,0}=\langle W\rangle.$ Given any
element $\Phi= a^1_1\oomega\wedge T+ a_1\oomega\wedge W+
a^1\orho\wedge T+ a\orho\wedge W$ in $\lie{g}^{*(0,1)}\otimes
\lie{g}^{1,0}$, $\dbar\Phi=\dbar(a^1\orho\wedge T)=
\frac{i}2a^1\oomega\wedge\orho\wedge W.$ Therefore, $
\ker\dbar\cap (\lie{g}^{*(0,1)}\otimes \lie{g}^{1,0})$ is the
linear span of $\oomega\wedge T$, $\oomega\wedge W$, and
$\orho\wedge W.$ Since the image of $\dbar$ in this space is
spanned by $\oomega\wedge W$, $\lie{h}^{1,1}=\langle \oomega\wedge
T, \orho\wedge W\rangle$. Similarly, kernel of $\dbar$ on the
space $\lie{g}^{*(0,2)}\otimes\lie{g}^{1,0}$ is spanned by
$\oomega\wedge\orho\wedge T$ and $\oomega\wedge\orho\wedge W$.
Since the latter is in the image of the $\dbar$-operator,
$\lie{h}^{1,2}=\langle \oomega\wedge\orho\wedge T\rangle.$

For the resolution of $\lie{g}^{2,0}$, since $\dbar(T\wedge W)=0$,
$\dbar\orho=0$ and $\dbar\oomega=0$,
\[
\lie{h}^{2,0} = \langle T\wedge W\rangle,\quad \lie{h}^{2,1}
=\langle \oomega\wedge T\wedge W,
    \orho\wedge T\wedge W\rangle,\quad
\lie{h}^{2,2} = \langle \oomega\wedge\orho\wedge T\wedge W
 \rangle.
\]
Due to Lemma \ref{full rank}, center of the algebra $\lie{h}$ with
respect to the Schouten bracket is the ideal generated by $\langle
\oomega\rangle$ with respect to the exterior product. In summary,
the space of harmonic representatives of the cohomology space
$\lie{h}$ is spanned by the following elements.

\

\begin{equation}\label{self brackets}
\begin{array}{|l|c|l|l|}\hline
 h^{p,q}&   &\mbox{in the center} & \mbox{in the complement of the center} \\\hline
 h^{0,0}& & 1 & \   \\\hline
 h^{0,1}& &\oomega &  \orho  \\\hline
 h^{0,2}&  & \oomega\wedge\orho & \ \\\hline
 h^{1,0}&  & W & \ \\\hline
 h^{1,1}&  & \oomega\wedge T & \orho\wedge W \\\hline
 h^{1,2}&  & \ & \orho\wedge T\wedge W \\\hline
 h^{2,0}&  & \ & T\wedge W  \\\hline
 h^{2,1}& & \oomega\wedge T\wedge W, \quad \oomega\wedge\orho\wedge T & \ \\\hline
 h^{2,2}&  & \oomega\wedge\orho\wedge T\wedge W & \ \\\hline
\end{array}
\end{equation}

\

By taking the above classes as harmonic representatives, we are in
effect choosing an inclusion of the Gerstenhaber algebra $\lie{h}$
in the full differential Gerstenhaber algebra $\lie{f}$. In this
sense, we consider orthogonal complement of $\lie{h}^{p,q}$ in
$\lie{f}^{p,q}$, denoted by $(\lie{h}^{p,q})^\perp$. This space is
spanned by the following elements.

\begin{equation}\label{outside h}
\begin{array}{|l|c|l|l|}\hline
 (h^{p,q})^\perp&   &\mbox{in the center} & \mbox{in the complement of the center} \\\hline
 f^{1,0}&  & \ & T \\\hline
 f^{1,1}&  & \oomega\wedge W=2i\dbar T & \orho\wedge T \\\hline
 f^{1,2}&  & \oomega\wedge \orho\wedge W=2i\dbar(\orho\wedge T) & \ \\\hline
\end{array}
\end{equation}

\

In the above computation, we repeatedly use the commutative law
(\ref{commutative}), distributive laws (\ref{distributive1}) and
(\ref{distributive2}), and
 information on degree-one elements to generate all the algebraic
relations among elements of higher degrees. It will be useful to
summarize all the \lq generating data\rq \ for the algebra.
\begin{lemma}\label{complex dG} The
DGA $dG(\lie{g}, J)$ $=(\lie{f}, \bra{}{}, \wedge, \dbar)$
associated to the complex structure of the Kodaira surface $(N,
J)$ is generated as follows.
\begin{enumerate}
\item The algebra with respect to the wedge product \lq$\wedge$' \
is the exterior algebra generated by  the degree-one elements $\{
T, W, \orho, \oomega\}$.
 \item The sole non-zero Schouten bracket
among degree-one elements is given by
$\bra{T}{\orho}=-\frac{i}2\oomega.$
 \item The differential on
degree-one elements is determined by $\dbar
T=-\frac{i}2\oomega\wedge W$,  $\dbar W=0$, $\dbar \orho=0$, and
$\dbar\oomega=0$.
\end{enumerate}
\end{lemma}

\begin{proposition}\label{abelian coho} The odd Lie superalgebra $(\lie{h}, \brah{}{})$ is Abelian.
\end{proposition}
\bproof Note that $\dbar W=0$. Therefore $W$ commutes with all
elements. Moreover for any $a$ and $b$ in $\lie{h}$, $
\bra{a}{b\wedge W}=\bra{a}{b}\wedge W.$ Therefore, the Schouten
brackets among the three elements $\orho\wedge W$, $\orho\wedge
T\wedge W$, and $T\wedge W$ are all equal to zero. We only need to
consider the Schouten bracket between $\orho$ and anyone of these
three elements. Due to Lemma \ref{2-dbar} and Lemma \ref{t-rho},
 the only non-zero brackets among elements in $\lie{h}$ are.
\begin{equation} \bra{\orho}{T\wedge W} =-\dbar
T \quad {\mbox{ and }} \quad \bra{\orho}{\orho\wedge T\wedge W}
=-\dbar (T\wedge\orho). \label{rho-rho-t-w}
\end{equation}
As they are $\dbar$-exact, the proof is completed. \eproof

Finally we complete the picture by presenting all the Schouten
brackets on $\lie{f}$ by including $T$ and $\orho\wedge T$. First
note that
\begin{equation}\label{o-n-o}
\bra{T}{\orho\wedge T}=-\frac{i}2\oomega\wedge T.
\end{equation}
The rest of the brackets is presented in the next table.

\

\begin{equation}\label{in-n-out}
\begin{array}{|l|c|l|l|}\hline
 \bra{A}{B}&   & T & \orho\wedge T \\\hline
\orho &  & \frac{i}2\oomega & -\frac{i}2\oomega\wedge\orho
\\\hline
\orho\wedge W & &\frac{i}2\oomega\wedge W=-\dbar T &
\frac{i}2\oomega\wedge\orho\wedge W=-\dbar (\orho\wedge T)
\\\hline
T\wedge W & &0 &\frac{i}2\oomega\wedge T\wedge W \\\hline
\orho\wedge T\wedge W& &\frac{i}2\oomega\wedge T\wedge W &
-i\oomega\wedge\orho\wedge T\wedge W
\\\hline
\end{array}
\end{equation}

\

\subsection{Solving Maurer-Cartan equation}\label{solving MC}
Consider an element $\Gamma_1$ in $\lie{h}$. We choose coordinate
functions as follows.
\begin{eqnarray}
 & &t_0+t_1\oomega\wedge\orho+t_2\orho\wedge
W+t_3\oomega\wedge T+t_4 T\wedge W +t_5\oomega\wedge\orho\wedge
T\wedge W\nonumber\\
&+& s_0\orho+s_1\oomega+s_2 W+s_3\orho\wedge T\wedge
W+s_4\oomega\wedge T\wedge W+s_5\oomega\wedge\orho\wedge
T.\label{coordinates}
\end{eqnarray}
Here $(t_0, \dots, t_5)$ are even coordinates and $(s_0, \dots,
s_5)$ are odd coordinates.

 In view of Table (\ref{in-n-out}) and Formula
(\ref{o-n-o}), we conclude that the only new terms could be added
through the recursive process are $T$ and $\orho \wedge T$.
Therefore, if $\Gamma$ is generated by the Kuranishi recursive
formula and $\Gamma_1$, there exist functions $\mu_1$ and $\mu_2$
such that
\[
\Gamma=\Gamma_1+\mu_1T+\mu_2\orho\wedge T.
\]
Next, we combine (\ref{rho-rho-t-w}), (\ref{o-n-o}) and
(\ref{in-n-out}) to conclude that
\begin{eqnarray}
& &\dbar\Gamma+\frac12\bra{\Gamma}{\Gamma}\nonumber\\
&=&\dbar\left((\mu_1-\mu_1t_2+s_0t_4)T\right)+\dbar\left((\mu_2-\mu_2t_2+s_0s_3)\orho\wedge
T\right)\nonumber\\
& &-\frac{i}2\mu_1s_0\oomega-\frac{i}2\mu_2s_0\oomega\wedge\orho
-i\mu_2s_3\oomega\wedge\orho\wedge T\wedge W-\frac{i}2\mu_1\mu_2\oomega\wedge T\nonumber\\
& &+\frac{i}2(\mu_2t_4-\mu_1s_3 )\oomega\wedge T\wedge W.
\label{temp}
\end{eqnarray}
Since $\oomega, \oomega\wedge\orho, \oomega\wedge\orho\wedge
T\wedge W$ and $\oomega\wedge T\wedge W$ are harmonic, $\Gamma$ is
a solution of the Maurer-Cartan equation only if
\[
\dbar\left((\mu_1-\mu_1t_2+s_0t_4)T\right)+\dbar\left((\mu_2-\mu_2t_2+s_0s_3)\orho\wedge
T\right)=0.
\]
When $|t_2|< 1$,  the solution is a pair of rational functions.
\[
\mu_1 =-\frac{s_0t_4}{1-t_2} \quad \mbox{ and } \quad  \mu_2
=-\frac{s_0s_3}{1-t_2}.
\]
In particular, $\mu_2t_4-\mu_1s_3=0.$ Therefore, when $|t_2|<1$,
\begin{equation}\label{gamma}
\Gamma=\Gamma_1-\frac{s_0t_4}{1-t_2}T-\frac{s_0s_3}{1-t_2}\orho\wedge
T.
\end{equation}
To complete a solution of the extended Maurer-Cartan equation, we
note that
\begin{eqnarray*}
\frac{\partial\Gamma}{\partial s_1}=\frac{\partial
\Gamma_1}{\partial s_1}=\oomega, & &\frac{\partial\Gamma}{\partial
t_1}=\frac{\partial
\Gamma_1}{\partial t_1}= \oomega\wedge\orho,  \\
\frac{\partial\Gamma}{\partial t_3}=\frac{\partial
\Gamma_1}{\partial t_3}=\oomega\wedge T, & &
\frac{\partial\Gamma}{\partial t_5}=\frac{\partial
\Gamma_1}{\partial t_5}= \oomega\wedge\orho\wedge T\wedge W.
\end{eqnarray*}
In view of (\ref{temp}), we set
\begin{eqnarray}
\vec{\partial}&=&\frac{i}2\mu_1s_0\frac{\partial}{\partial s_1}
+\frac{i}2\mu_2s_0\frac{\partial}{\partial
t_1}+\frac{i}2\mu_1\mu_2\frac{\partial}{\partial t_3}
 +i\mu_2s_3\frac{\partial}{\partial t_5}\nonumber\\
 &=&-\frac{i}2\frac{s_0}{1-t_2}
 \left(t_4s_0\frac{\partial}{\partial s_1}
 +s_3s_0\frac{\partial}{\partial t_1}
 -\frac{t_4s_0s_3}{1-t_2}\frac{\partial}{\partial t_3}
 +2s^2_3\frac{\partial}{\partial t_5}\right).\label{chen}
 \end{eqnarray}
 We sum up the above observations below.
 \begin{theorem}\label{general solution} Suppose that $\Gamma_1$ is an element in the Dolbeault
 cohomology of polyvector fields on the Kodaira surface. Choose a
 harmonic basis such that $\Gamma_1$ is given by
 {\rm (\ref{coordinates})}. Then $\Gamma$ given by {\rm (\ref{gamma})} and
 $\vec\partial$ given by {\rm (\ref{chen})}  solve
 the extended Maurer-Cartan equation {\rm (\ref{EMC})}.
 \end{theorem}

\subsection{Extended Kuranishi space}\label{extended k}
Due to (\ref{chen}), the Kuranishi space for Kodaira surface has
two components, namely
\begin{equation}
{\ck_0}=\{s_0=0\} \quad \mbox{ and } \quad {\ck_1}=\{t_4=0, \quad
s_3=0\}.
\end{equation}

When $s_0=0$, $\Gamma=\Gamma_1$.  Since all the coefficient
functions for the coordinate vector fields except
$\frac{\partial}{\partial t_5}$ in $\vec\partial$ vanish along
$s_0=0$ to order 2, $[\vec\partial, Y]$
 is non-zero along $s_0=0$ only if $Y=\frac{\partial}{\partial
s_0}$. In this case,
\begin{equation}
[\vec\partial, \frac{\partial}{\partial s_0}]
=-i\frac{s_3^2}{1-t_2}\frac{\partial}{\partial t_5} \quad \mbox{
mod } s_0.
\end{equation}
It follows that the Kuranishi space $\ck_0$ has two strata. One is
a linear space defined by $s_3=0$. It is isomorphic to
$\CC^{6|4}$. The other stratum is the quotient of the open set
$s_3\neq 0$ with respect to the distribution spanned by
$\frac{\partial}{\partial t_5}$. Therefore,
 it is identified to the super-space
\begin{equation}
\ck_{0, \mbox{generic}}=\{(t_0, t_1, t_2, t_3, t_4, s_1, s_2, s_3,
s_4, s_5)\in \CC^{5|5}: s_3\neq 0\}.
\end{equation}

The restriction of the distribution $\cd$ to $\ck_1$ is spanned by
\[
2i[\vec\partial, \frac{\partial}{\partial s_3}]=
-\frac{s^2_0}{1-t_2}\frac{\partial}{\partial t_1} \quad \mbox{ and
} \quad 2i[\vec\partial, \frac{\partial}{\partial t_4}]=
\frac{s^2_0}{1-t_2}\frac{\partial}{\partial s_1}.
\]
It has two strata. A stratum is the linear subspace of $\lie{h}$
defined by $s_0=0$. It is isomorphic to $\CC^{5|4}$.
 The component given by $s_0\neq
0$ is the quotient of an open set in $\CC^{5|5}$ with respect to
$\frac{\partial}{\partial t_1}$ and $\frac{\partial}{\partial
s_1}$. Therefore, it is an open set contained in $\CC^{4|4}$.

\begin{theorem}\label{kur components}
 The extended Kuranishi space of the
Kodaira surface has four components. There are  two linear
components $\CC^{6|4}$ and $\CC^{5|4}$. It has two additional
components contained in the complement of an odd linear hyperplane
in $\CC^{5|5}$ and in $\CC^{4|4}$ respectively.
\end{theorem}

\subsection{Associative product for a weak Frobenius manifold}\label{find mu}
 In this section, we  calculate the associative product
on the weak Frobenius manifold $\ck_{0, \mbox{generic}}$ as
explained in Section \ref{review}. We choose this component to
compute partly because of the simplicity due to the identity
$\Gamma=\Gamma_1$.

Since $\frac{\partial \Gamma}{\partial t_0}=1$, the tangent vector
field $\frac{\partial }{\partial t_0}$ is the unit element. Next
we calculate the wedge products among
\[
\frac{\partial \Gamma}{\partial t_\alpha}, \quad \frac{\partial
\Gamma}{\partial s_\beta}, \quad \alpha=1,2,3,4, \quad
\beta=1,2,3,4,5.
\]
For example,
\[
\frac{\partial \Gamma}{\partial t_1}\wedge \frac{\partial
\Gamma}{\partial t_4} =(\oomega\wedge\orho)\wedge (T\wedge W)
=\frac{\partial \Gamma}{\partial t_5}. \] Since we take quotient
respect to $\frac{\partial \Gamma}{\partial t_5}$, $
\frac{\partial }{\partial t_1}\circ \frac{\partial }{\partial t_4}
=0. $ Let us take two more examples. Consider
\[
\frac{\partial \Gamma}{\partial s_1}\wedge \frac{\partial
\Gamma}{\partial s_2}=\oomega\wedge W. \]
Since
\[
\dbar_\Gamma(T) =-\frac{i}2\oomega\wedge
W+\frac{i}2t_2\oomega\wedge W +\frac{i}2s_3\oomega\wedge T\wedge
W,
\]
\[
(1-t_2)\oomega\wedge W\equiv s_3\oomega\wedge T\wedge W ={s_3}
\frac{\partial \Gamma}{\partial s_4} \quad \mbox{ mod }
\image\dbar_\Gamma. \]
Therefore,
\[
\frac{\partial \Gamma}{\partial s_1}\wedge \frac{\partial
\Gamma}{\partial s_2}\equiv \frac{s_3}{1-t_2} \frac{\partial
\Gamma}{\partial s_4} \quad \mbox{ mod } \image\dbar_\Gamma, \]
and
\[
\frac{\partial }{\partial s_1}\circ \frac{\partial}{\partial
s_2}=\frac{s_3}{1-t_2} \frac{\partial }{\partial s_4}.
\]

Next, we consider
\[
\frac{\partial \Gamma}{\partial t_1}\wedge \frac{\partial
\Gamma}{\partial s_2}=\oomega\wedge\orho\wedge W. \]
Given
$s_0=0$,
\begin{eqnarray*}
\dbar_\Gamma(\orho\wedge T)&=&-\frac{i}2\oomega\wedge\orho\wedge
W+\frac{i}2t_2\oomega\wedge \orho\wedge W
+\frac{i}2t_4\oomega\wedge T\wedge W
\nonumber\\
& &-is_3\oomega\wedge\orho \wedge T\wedge W.
\end{eqnarray*}
Therefore,
\begin{eqnarray*}
& &(1-t_2)\oomega\wedge\orho\wedge W\\
&\equiv& t_4\oomega\wedge T\wedge W -2s_3\oomega\wedge\orho \wedge
T\wedge W =t_4\frac{\partial \Gamma}{\partial
s_4}-2s_3\frac{\partial\Gamma}{\partial t_5} \quad \mbox{ mod }
\image\dbar_\Gamma.
\end{eqnarray*}
Since we take quotient with respect to $\frac{\partial}{\partial
t_5}$,
\[
\frac{\partial \Gamma}{\partial t_1}\wedge \frac{\partial
\Gamma}{\partial s_2}\equiv\frac{t_4}{1-t_2}\frac{\partial
\Gamma}{\partial s_4} \quad \mbox{ mod } \image \dbar_\Gamma,
\frac{\partial}{\partial t_5}. \]
Therefore,
\[
\frac{\partial}{\partial t_1}\circ \frac{\partial }{\partial
s_2}=\frac{t_4}{1-t_2}\frac{\partial }{\partial s_4}. \]
We apply
the above computation to every pair of coordinate vectors on the
Kuranishi space $\ck_{0, \mbox{generic}}$. It turns out that
$\frac{\partial}{\partial s_3}, \frac{\partial}{\partial s_4},
\frac{\partial}{\partial s_5}$ do not have non-trivial product.
Recall that $\frac{\partial}{\partial t_0}$ is the unit element
with respect to $\circ$. We put non-trivial product in a table
below with respect to the ordered bases
\[
\frac{\partial}{\partial t_1}, \frac{\partial}{\partial t_2},
\frac{\partial}{\partial t_3}, \frac{\partial}{\partial t_4},
\frac{\partial}{\partial s_1}, \frac{\partial}{\partial s_2}.
\]
\begin{equation}
\begin{array}{|c|c|c|c|c|c|c|c|}\hline
 \circ& \frac{\partial}{\partial t_1} & \frac{\partial}{\partial
 t_2}&
\frac{\partial}{\partial t_3} & \frac{\partial}{\partial t_4} &
\frac{\partial}{\partial s_1} & \frac{\partial}{\partial s_2}
  \\\hline
  \frac{\partial}{\partial t_1} & 0& 0&0&0 & 0& \frac{t_4}{1-t_2}\frac{\partial}{\partial s_4}
\\\hline
   \frac{\partial}{\partial t_2}& 0& 0&0&0 &\frac{t_4}{1-t_2}\frac{\partial}{\partial
   s_4}&0
\\\hline
\frac{\partial}{\partial t_3} & 0& 0&0&0 &
0&\frac{\partial}{\partial s_4}\\\hline \frac{\partial}{\partial
t_4}& 0& 0&0&0 &\frac{\partial}{\partial s_4}& 0\\\hline
 \frac{\partial}{\partial s_1}& 0&\frac{t_4}{1-t_2}\frac{\partial}{\partial
 s_4}&0&\frac{\partial}{\partial s_4}&0& \frac{s_3}{1-t_2}\frac{\partial}{\partial
 s_4}\\\hline
\frac{\partial}{\partial s_2}
&\frac{t_4}{1-t_2}\frac{\partial}{\partial
 s_4}&0&\frac{\partial}{\partial
 s_4}&0&-\frac{s_3}{1-t_2}\frac{\partial}{\partial
 s_4}&0
\\\hline
\end{array}
\end{equation}

\section{Mirror Image of Kodaira Surfaces}
 In this
section, we continue our analysis of extended deformation of a
complex structure on the Kodaira surface, but limit to harmonic
2-fields as initial conditions.

\subsection{From complex to symplectic structures}\label{to
symplectic}  Recall that $\lie{h}^2$ is spanned by
$\cb=\oomega\wedge\orho$, $\psi=\orho\wedge W$,
$\phi=\oomega\wedge T$, and $B=T\wedge W$. If $\Gamma$ is
contained in $\lie{h}^2$, there are complex numbers $(t_1, t_2,
t_3, t_4)$ such that
\begin{equation}\label{reduced}
\Gamma=t_1\cb+t_2\psi+t_3\phi+t_4B.
\end{equation}
Due to dimension limitation, the non-harmonic classical field
$s^i_j$ vanishes. Therefore, by either Theorem \ref{generalized}
or Theorem  \ref{general solution}, every element in $\lie{h}^2$
together with $\overrightarrow{\partial}=0$ is a solution of the
Maurer-Cartan equation.

We consider the following interpretation of the 2-fields.
\begin{eqnarray*}
\cb &\in& \lie{g}^{*(0,2)}\subset \End (\lie{g}^{0,1}, \lie{g}^{*(0,1)}).\\
\psi, \phi &\in& \lie{g}^{*(0,1)}\otimes \lie{g}^{1,0}\subset \End(\lie{g}^{(0,1)}, \lie{g}^{1,0})
\equiv\End (\lie{g}^{*(1,0)},\lie{g}^{*(0,1)}).\\
B &\in& \lie{g}^{2,0}\subset \End(\lie{g}^{*(1,0)},
\lie{g}^{1,0}).
\end{eqnarray*}
In such context, we see that
\begin{eqnarray}
\cb(\overline{T}) &=& \orho, \quad \cb(\overline{W}) =-\oomega,
\quad B(\omega)=W, \quad
B(\rho)=-T\\
\psi(\overline{W}) &=&W, \quad \psi(\rho)=-\orho, \quad
\phi(\overline{T})=T, \quad \phi(\omega)=-\oomega.
\end{eqnarray}
Now consider  the distribution spanned by ${\overline L}=\langle
\overline{T}, \overline{W}, \omega, \rho\rangle.$ It is considered
as a subbundle of $(T_N\oplus T_N^*)_\CC$. This is a choice of a
complex structure through the distribution of (0,1)-vectors and
their annihilators. It defines  a generalized complex structure
\cite{Gua}.  Given $\Gamma$ as in (\ref{reduced}) a new
distribution ${\overline L}_{\Gamma}$ is defined by its graph.
With respect to the ordered base $ \{\overline{T}, \overline{W},
\omega, \rho, T, W, \oomega, \orho\}, $ the distribution
${\overline{L}}_\Gamma$ is given by the
\begin{equation}\label{frame}
\left(
\begin{array}
[c]{c}%
\overline{T}+\Gamma(\overline{T})\\
\overline{W}+\Gamma(\overline{W})\\
\omega+\Gamma(\omega)\\
\rho+\Gamma(\rho)
\end{array}
\right)  =\left(
\begin{array}
[c]{cccccccc}%
1 & 0 & 0 & 0 & t_{3} & 0 & 0 & t_{1}\\
0 & 1 & 0 & 0 & 0 & t_{2} & -t_{1} & 0\\
0 & 0 & 1 & 0 & 0 & t_{4} & -t_{3} & 0\\
0 & 0 & 0 & 1 & -t_{4} & 0 & 0 & -t_{2}%
\end{array}
\right) .
\end{equation}

If $t_1=t_4=0$, we recover classical deformation of complex
structures on Kodaira surfaces as found by Borcea \cite{Borcea}
\cite{GMPP}. We may also consider the extended complex structure
given by $t_2=t_3=0$.  There is a very special family depending on
one complex parameter.
\begin{equation}
{\overline L}_\Gamma=\langle \overline{T}+\frac12 t\orho, \quad
\overline{W}-\frac12 t\oomega, \quad \omega-\frac{2}{\overline
t}W, \quad \rho+\frac{2}{\overline t}T \rangle,
\end{equation}
where $t\neq 0$. Let $u$ and $v$ be the real and imaginary part of
the complex number $t$. Let $\omega=\alpha+i\beta$ and
$\rho=\gamma+i\delta$. Define
\begin{equation}
\Omega=u(\alpha\wedge\gamma-\beta\wedge\delta)
  +v(\alpha\wedge\delta+\beta\wedge\gamma).
\end{equation}
Then the distribution ${\overline L}_\Gamma\oplus L_\Gamma$ is the
complexification of
\begin{equation}
\langle X+\iota_X\Omega, \quad Y+\iota_Y\Omega, \quad
U+\iota_U\Omega, \quad V+\iota_V\Omega \rangle
\end{equation}
as a subbundle of the real direct sum $T_N\oplus T_N^*$.
Therefore, the complex structure $J$ deforms through extended
deformation theory to an invariant symplectic structure $\Omega$
in the form of a generalized complex structure \cite{Gua}. This
family of symplectic structures is contained in the family of all
invariant symplectic structures on $N$ as we shall see next.

\subsection{DGA of symplectic
structures}\label{def sym} For any set of real numbers $(u_1, v_1,
u_2, v_2)$ with $\triangle:=u_1^2+v_1^2-u_2^2-v_2^2\neq 0$, we
define $\Omega$ to be the closed 2-form below.
\begin{equation}\label{general omega}
u_1(\alpha\wedge\gamma-\beta\wedge\delta)+v_1(\alpha\wedge\delta+\beta\wedge\gamma)
+u_2(\alpha\wedge\gamma+\beta\wedge\delta)+v_2(\alpha\wedge\delta-\beta\wedge\gamma),
\end{equation}
Since $\triangle\neq 0$, $\Omega$ is a symplectic form on the
nilmanifold $N$.

Given the symplectic structure, we obtain a DGA in the standard
way. Namely, a  contraction with the symplectic form $\Omega$
defines an isomorphism from the tangent bundle to the cotangent
bundle on the nilmanifold $N$. Then the Lie bracket among vector
fields is carried by the inverse isomorphism to a Schouten bracket
$\bra{}{}_\Omega$ on de Rham algebra of differential forms on $N$.
The package
\[
dG(N, \Omega)=\left(\oplus_\ell C^\infty(N, \wedge^\ell T^*_N),
\bra{}{}_\Omega, \wedge, d \right)
\]
is a DGA   \cite{Mer-semi}. Its cohomology is the de Rham
cohomology. Such construction could be limited to the space of
invariant differential forms. It yields a new DGA $dG(\lie{g},
\Omega)$ $=\left( \wedge^*\lie{g}^*, \bra{}{}_\Omega, \wedge, d
\right)$. Due to Nomizu, the inclusion of the complex of invariant
differential forms of nilmanifolds in the de Rham complex induces
an isomorphism of cohomology \cite{Nomizu}. In other words, we
have
\begin{proposition}\label{symp qua}
The inclusion of invariant algebra $dG(\lie{g}, \Omega)$  in the
algebra $dG(N, \Omega)$ is a quasi-isomorphism.
\end{proposition}
In our case, the contraction $\iota$ with $\Omega$ yields the
following.
\begin{eqnarray*}
\iota(X)&=&(u_1+u_2)\gamma+(v_1+v_2)\delta, \quad \iota(Y)=(v_1-v_2)\gamma-(u_1-u_2)\delta,\\
\iota(U)&=&-(u_1+u_2)\alpha-(v_1-v_2)\beta, \quad
\iota(V)=-(v_1+v_2)\alpha+(u_1-u_2)\beta.
\end{eqnarray*}
Let $\alpha'=-\frac{1}{\triangle}\iota(U)$ and
$\beta'=\frac{1}{\triangle}\iota(V)$. It follows that $
d\gamma=-\alpha\wedge\beta=-\triangle \alpha'\wedge\beta'.$ As the
only non-trivial bracket among the vectors $X,Y,U,V$ is $[X,Y]=U$,
the only non-zero bracket among the 1-forms
$\alpha',\beta',\gamma,\delta$ is given by
\begin{equation}
\bra{\gamma}{\delta}_\Omega=\alpha'.
\end{equation}
\begin{lemma}\label{symp dG}
Let $\Omega$ be the symplectic form on the Kodaira surface given
by {\rm (\ref{general omega})}. The invariant differential
Gerstenhaber algebra $dG(\lie{g}, \Omega)$ is generated as
follows.
\begin{enumerate}
\item The algebra with respect to the wedge product is the
exterior algebra generated by  the degree-one elements $\alpha',
\beta', \gamma, \delta$. \item The sole non-zero Schouten bracket
among degree-one elements is given by
$\bra{\gamma}{\delta}_\Omega=\alpha'.$
 \item The differential on
degree-one elements is determined by $d\gamma=-\triangle
\alpha'\wedge \beta',$ $d\beta'=0,$ $d\delta=0,$ and $d\alpha'=0.$
\end{enumerate}
\end{lemma}

\begin{theorem}\label{follow-up} Let $N$ be the Kodaira surface. There is  an isomorphism
\[
\Upsilon: (\oplus_{k=0}^2\sum_{p+q=k}H^q(N, \wedge^p\ct_N),
\bra{}{}, \wedge ) \to \left(\oplus_{k=0}^2H^k(N, \CC),
\bra{}{}_\Omega, \wedge \right).
\]
\end{theorem}
\bproof In view  of Lemma \ref{complex dG} and Lemma \ref{symp
dG},  the map $\Upsilon$ defined by $\dbar\mapsto d,$ $T\mapsto
\gamma$, $W\mapsto -\triangle\beta'$, $\orho \mapsto \delta$,
$-\frac{i}2\oomega \mapsto \alpha'$ produces an isomorphism from
$dG(\lie{g}, J)$ to $dG(\lie{g}, \Omega)$. The theorem then
follows the quasi-isomorphisms established by Proposition \ref{cx
qua} and Proposition \ref{symp qua}. \eproof

\noindent{\bf Acknowledgment\ } \  I thank Sergei Merkulov for
explaining his work to me. I thank Paul Gauduchon directing my
attention to Gualtieri's thesis. Lots of techniques concerned with
cohomology on nilmanifolds are due to my past collaboration with
Gueo Grantcharov and Henrik Pedersen, and subsequently with Colin
Maclaughlin and Simon Salamon. I also thank J. Zhou and N. C.
Leung for useful comments.

\vspace{.25in}

\noindent\texttt{ypoon@ucr.edu}\\
Corresponding Address: Department of Mathematics, University of
California at Riverside, Riverside, CA 92521, USA.

\end{document}